\documentclass[11pt]{amsart}
\usepackage{amsmath}
\usepackage{amssymb,latexsym}
\usepackage{amscd}
\usepackage{xspace}
\usepackage{verbatim}
\begin{document}
\numberwithin{equation}{section}
\title[Classification of positive solutions]{Complete classification of the positive solutions of
$-\Gd u+ u^q=0$}
\author{Moshe Marcus }
\address{Department of Mathematics, Technion\\
 Haifa 32000, ISRAEL}
 \email{marcusm@math.technion.ac.il}

\date{\today}

\newcommand{\txt}[1]{\;\text{ #1 }\;}
\newcommand{\tbf}{\textbf}
\newcommand{\tit}{\textit}
\newcommand{\tsc}{\textsc}
\newcommand{\trm}{\textrm}
\newcommand{\mbf}{\mathbf}
\newcommand{\mrm}{\mathrm}
\newcommand{\bsym}{\boldsymbol}
\newcommand{\scs}{\scriptstyle}
\newcommand{\sss}{\scriptscriptstyle}
\newcommand{\txts}{\textstyle}
\newcommand{\dsps}{\displaystyle}
\newcommand{\fnz}{\footnotesize}
\newcommand{\scz}{\scriptsize}
\newcommand{\be}{\begin{equation}}
\newcommand{\bel}[1]{\begin{equation}\label{#1}}
\newcommand{\ee}{\end{equation}}
\newtheorem{subn}{\name}
\renewcommand{\thesubn}{}
\newcommand{\bsn}[1]{\def\name{#1$\!\!$}\begin{subn}}
\newcommand{\esn}{\end{subn}}
\newtheorem{sub}{\name}[section]
\newcommand{\dn}[1]{\def\name{#1}}   
\newcommand{\bs}{\begin{sub}}
\newcommand{\es}{\end{sub}}
\newcommand{\bsl}[1]{\begin{sub}\label{#1}}
\newcommand{\bth}[1]{\def\name{Theorem}\begin{sub}\label{t:#1}}
\newcommand{\blemma}[1]{\def\name{Lemma}\begin{sub}\label{l:#1}}
\newcommand{\bcor}[1]{\def\name{Corollary}\begin{sub}\label{c:#1}}
\newcommand{\bdef}[1]{\def\name{Definition}\begin{sub}\label{d:#1}}
\newcommand{\bprop}[1]{\def\name{Proposition}\begin{sub}\label{p:#1}}
\newcommand{\bnote}[1]{\def\name{\mdseries\scshape Notation}\begin{sub}\label{n:#1}}
\newcommand{\bproof}{\begin{proof}}
\newcommand{\eproof}{\end{proof}}
\newcommand{\bcom}{}
\newcommand{\req}{\eqref}
\newcommand{\rth}[1]{Theorem~\ref{t:#1}}
\newcommand{\rlemma}[1]{Lemma~\ref{l:#1}}
\newcommand{\rcor}[1]{Corollary~\ref{c:#1}}
\newcommand{\rdef}[1]{Definition~\ref{d:#1}}
\newcommand{\rprop}[1]{Proposition~\ref{p:#1}}
\newcommand{\rnote}[1]{Notation~\ref{n:#1}}
\newcommand{\BA}{\begin{array}}
\newcommand{\EA}{\end{array}}
\newcommand{\BAN}{\renewcommand{\arraystretch}{1.2}
\setlength{\arraycolsep}{2pt}\begin{array}}
\newcommand{\BAV}[2]{\renewcommand{\arraystretch}{#1}
\setlength{\arraycolsep}{#2}\begin{array}}
\newcommand{\BSA}{\begin{subarray}}
\newcommand{\ESA}{\end{subarray}}
\newcommand{\BAL}{\begin{aligned}}
\newcommand{\EAL}{\end{aligned}}
\newcommand{\BALG}{\begin{alignat}}
\newcommand{\EALG}{\end{alignat}}
\newcommand{\BALGN}{\begin{alignat*}}
\newcommand{\EALGN}{\end{alignat*}}
\newcommand{\note}[1]{\noindent\textit{#1.}\hspace{2mm}}
\newcommand{\Remark}{\note{Remark}}

\newcommand{\forevery}{\quad \forall}
\newcommand{\1}{\\[1mm]}
\newcommand{\2}{\\[2mm]}
\newcommand{\3}{\\[3mm]}
\newcommand{\set}[1]{\{#1\}}
\def\({{\rm (}}
\def\){{\rm )}}
\newcommand{\st}[1]{{\rm (#1)}}
\newcommand{\lra}{\longrightarrow}
\newcommand{\lla}{\longleftarrow}
\newcommand{\llra}{\longleftrightarrow}
\newcommand{\Lra}{\;\;\Longrightarrow\;\;}
\newcommand{\Lla}{\;\;\Longleftarrow\;\;}
\newcommand{\Llra}{\Longleftrightarrow}
\newcommand{\harpup}{\rightharpoonup}
\def\dar{\downarrow}
\def\uar{\uparrow}
\newcommand{\paran}[1]{\left (#1 \right )}
\newcommand{\sqrbr}[1]{\left [#1 \right ]}
\newcommand{\curlybr}[1]{\left \{#1 \right \}}
\newcommand{\abs}[1]{\left |#1\right |}
\newcommand{\norm}[1]{\left \|#1\right \|}
\newcommand{\angbr}[1]{\left< #1\right>}
\newcommand{\paranb}[1]{\big (#1 \big )}
\newcommand{\sqrbrb}[1]{\big [#1 \big ]}
\newcommand{\curlybrb}[1]{\big \{#1 \big \}}
\newcommand{\absb}[1]{\big |#1\big |}
\newcommand{\normb}[1]{\big \|#1\big \|}
\newcommand{\angbrb}[1]{\big\langle #1 \big \rangle}
\newcommand{\thkl}{\rule[-.5mm]{.3mm}{3mm}}
\newcommand{\thknorm}[1]{\thkl #1 \thkl\,}
\newcommand{\trinorm}[1]{|\!|\!| #1 |\!|\!|\,}
\newcommand{\vstrut}[1]{\rule{0mm}{#1}}
\newcommand{\rec}[1]{\frac{1}{#1}}
\newcommand{\opname}[1]{\mathrm{#1}\,}
\newcommand{\supp}{\opname{supp}}
\newcommand{\dist}{\opname{dist}}
\newcommand{\sign}{\opname{sign}}
\newcommand{\diam}{\opname{diam}}
\newcommand{\proj}{\opname{proj}}
\newcommand{\q}{\quad}
\newcommand{\qq}{\qquad}
\newcommand{\hsp}[1]{\hspace{#1mm}}
\newcommand{\vsp}[1]{\vspace{#1mm}}
\newcommand{\prt}{\partial}
\newcommand{\sms}{\setminus}
\newcommand{\ems}{\emptyset}
\newcommand{\ti}{\times}
\newcommand{\pr}{^\prime}
\newcommand{\ppr}{^{\prime\prime}}
\newcommand{\tl}{\tilde}
\newcommand{\wtl}{\widetilde}
\newcommand{\sbs}{\subset}
\newcommand{\sbeq}{\subseteq}
\newcommand{\nind}{\noindent}
\newcommand{\ovl}{\overline}
\newcommand{\unl}{\underline}
\newcommand{\nin}{\not\in}
\newcommand{\pfrac}[2]{\genfrac{(}{)}{}{}{#1}{#2}}
\newcommand{\tin}{\to\infty}
\newcommand{\ind}[1]{_{_{#1}}\!}
\newcommand{\chr}[1]{\chi\ind{#1}}
\newcommand{\rest}[1]{\big |\ind{#1}}
\newcommand{\num}[1]{{\rm (#1)}\hspace{2mm}}
\newcommand{\wkc}{weak convergence\xspace}
\newcommand{\wrto}{with respect to\xspace}
\newcommand{\cons}{consequence\xspace}
\newcommand{\consy}{consequently\xspace}
\newcommand{\Consy}{Consequently\xspace}
\newcommand{\Essy}{Essentially\xspace}
\newcommand{\essy}{essentially\xspace}
\newcommand{\mnz}{minimizer\xspace}
\newcommand{\sth}{such that\xspace}
\newcommand{\ngh}{neighborhood\xspace}
\newcommand{\nghs}{neighborhoods\xspace}
\newcommand{\seq}{sequence\xspace}
\newcommand{\sseq}{subsequence\xspace}
\newcommand{\locun}{locally uniformly\xspace}
\newcommand{\ifif}{if and only if\xspace}
\newcommand{\suff}{sufficiently\xspace}
\newcommand{\abc}{absolutely continuous\xspace}
\newcommand{\sol}{solution\xspace}
\newcommand{\subsol}{subsolution\xspace}
\newcommand{\supsol}{supersolution\xspace}
\newcommand{\Wlg}{Without loss of generality\xspace}
\newcommand{\wlg}{without loss of generality\xspace}
\newcommand{\bdw}{\partial\Gw}
\newcommand{\Capq}{C_{2/q,q'}}
\newcommand{\Cqq}{$C_{2/q,q'}$}
\newcommand{\finecl}{$C_{2/q,q'}$-finely closed \xspace}
\newcommand{\fineop}{$C_{2/q,q'}$-finely open \xspace}
\newcommand{\finetop}{$C_{2/q,q'}$-fine topology \xspace}
\newcommand{\Cqconv}{$C_{2/q,q'}$-convergent \xspace}

\newcommand{\Wq}{W^{2/q,q'}}
\newcommand{\Wqq}{W^{-2/q,q}}
\newcommand{\Wqdb}{W^{-2/q,q}_+(\bdw)}
\newcommand{\sbsq}{\overset{q}{\sbs}}
\newcommand{\smq}{\overset{q}{\sim}}
\newcommand{\app}[1]{\underset{#1}{\approx}}
\newcommand{\suppq}{\mathrm{supp}^q_{\bdw}\,}
\newcommand{\convq}{\overset{q}{\to}}
\newcommand{\barq}[1]{\bar{#1}^{^q}}
\newcommand{\prtq}{\partial_q}
\newcommand{\tr}{\mathrm{tr}\,}
\newcommand{\Tr}{\mathrm{Tr}\,}
\newcommand{\trR}{\mathrm{tr}\ind{\CR}}
\newcommand{\trin}[1]{\mathrm{tr}\ind{#1}}
\newcommand{\qcl}{$\Capq$-finely closed\xspace}
\newcommand{\qop}{$\Capq$-finely open\xspace}
\newcommand{\gsmod}{$\gs$-moderate\xspace}
\newcommand{\gsreg}{$\gs$-regular\xspace}
\newcommand{\qreg}{$q$-quasi regular\xspace}
\newcommand{\qeq}{$\Capq$-equivalent\xspace}
\newcommand{\ppf}{\underset{f}{\prec\prec}}
\newcommand{\ofrown}{\overset{\frown}}
\newcommand{\modcon}{\underset{mod}{\lra}}
\newcommand{\ugb}[1]{u\chr{\Gs_\gb(#1)}}
\newcommand{\mcon}{$q$-moderately convergent\xspace}
\newcommand{\mdiv}{$q$-moderately divergent\xspace}
\def\qsupp{q\text{-supp}\,}
\def\Lim{\,\text{\rm Lim}\,}
\def\muCR{\mu\ind{\CR}}
\def\vCR{v\ind{\CR}}
\def\bcom{}
\def\ga{\alpha}     \def\gb{\beta}       \def\gg{\gamma}
\def\gc{\chi}       \def\gd{\delta}      \def\ge{\epsilon}
\def\gth{\theta}                         \def\vge{\varepsilon}
\def\gf{\phi}       \def\vgf{\varphi}    \def\gh{\eta}
\def\gi{\iota}      \def\gk{\kappa}      \def\gl{\lambda}
\def\gm{\mu}        \def\gn{\nu}         \def\gp{\pi}
\def\vgp{\varpi}    \def\gr{\rho}        \def\vgr{\varrho}
\def\gs{\sigma}     \def\vgs{\varsigma}  \def\gt{\tau}
\def\gu{\upsilon}   \def\gv{\vartheta}   \def\gw{\omega}
\def\gx{\xi}        \def\gy{\psi}        \def\gz{\zeta}
\def\Gg{\Gamma}     \def\Gd{\Delta}      \def\Gf{\Phi}
\def\Gth{\Theta}
\def\Gl{\Lambda}    \def\Gs{\Sigma}      \def\Gp{\Pi}
\def\Gw{\Omega}     \def\Gx{\Xi}         \def\Gy{\Psi}

\def\CS{{\mathcal S}}   \def\CM{{\mathcal M}}   \def\CN{{\mathcal N}}
\def\CR{{\mathcal R}}   \def\CO{{\mathcal O}}   \def\CP{{\mathcal P}}
\def\CA{{\mathcal A}}   \def\CB{{\mathcal B}}   \def\CC{{\mathcal C}}
\def\CD{{\mathcal D}}   \def\CE{{\mathcal E}}   \def\CF{{\mathcal F}}
\def\CG{{\mathcal G}}   \def\CH{{\mathcal H}}   \def\CI{{\mathcal I}}
\def\CJ{{\mathcal J}}   \def\CK{{\mathcal K}}   \def\CL{{\mathcal L}}
\def\CT{{\mathcal T}}   \def\CU{{\mathcal U}}   \def\CV{{\mathcal V}}
\def\CZ{{\mathcal Z}}   \def\CX{{\mathcal X}}   \def\CY{{\mathcal Y}}
\def\CW{{\mathcal W}}
\def\BBA {\mathbb A}   \def\BBb {\mathbb B}    \def\BBC {\mathbb C}
\def\BBD {\mathbb D}   \def\BBE {\mathbb E}    \def\BBF {\mathbb F}
\def\BBG {\mathbb G}   \def\BBH {\mathbb H}    \def\BBI {\mathbb I}
\def\BBJ {\mathbb J}   \def\BBK {\mathbb K}    \def\BBL {\mathbb L}
\def\BBM {\mathbb M}   \def\BBN {\mathbb N}    \def\BBO {\mathbb O}
\def\BBP {\mathbb P}   \def\BBR {\mathbb R}    \def\BBS {\mathbb S}
\def\BBT {\mathbb T}   \def\BBU {\mathbb U}    \def\BBV {\mathbb V}
\def\BBW {\mathbb W}   \def\BBX {\mathbb X}    \def\BBY {\mathbb Y}
\def\BBZ {\mathbb Z}

\def\GTA {\mathfrak A}   \def\GTB {\mathfrak B}    \def\GTC {\mathfrak C}
\def\GTD {\mathfrak D}   \def\GTE {\mathfrak E}    \def\GTF {\mathfrak F}
\def\GTG {\mathfrak G}   \def\GTH {\mathfrak H}    \def\GTI {\mathfrak I}
\def\GTJ {\mathfrak J}   \def\GTK {\mathfrak K}    \def\GTL {\mathfrak L}
\def\GTM {\mathfrak M}   \def\GTN {\mathfrak N}    \def\GTO {\mathfrak O}
\def\GTP {\mathfrak P}   \def\GTR {\mathfrak R}    \def\GTS {\mathfrak S}
\def\GTT {\mathfrak T}   \def\GTU {\mathfrak U}    \def\GTV {\mathfrak V}
\def\GTW {\mathfrak W}   \def\GTX {\mathfrak X}    \def\GTY {\mathfrak Y}
\def\GTZ {\mathfrak Z}   \def\GTQ {\mathfrak Q}
\font\Sym= msam10
\def\SYM#1{\hbox{\Sym #1}}

\def\bmn{\mathbf{n}}
\def\bma{\mathbf{a}}
\newcommand{\prtn}{\prt_{\mathbf{n}}}
\def\txin{\txt{in}}
\def\txon{\txt{on}}
\def\rhs{\text{right hand side \xspace}}
\def\lhs{\text{left hand side \xspace}}
\def\L1wsol{$L^1$-weak solution\xspace}
\def\W1p{W^{1,p}}
\def\Lr{L^1(\Gw,\gr)}
\def\lfs#1{\lfloor_{\sss  #1}}
\def\GTMr{\GTM(\Gw;\gr)}
\def\bgw{\bar\Gw}
\def\prtn{\prt_{\bmn}}
\def\loc{_{\mathrm{loc}}}
\def\bvp{boundary value problem\xspace}
\def\superh{superharmonic\xspace}
\def\subh{subharmonic\xspace}
\def\RN{\BBR^N}
\def\HMV{H^V_{\textnormal{mod}}}
\begin{abstract}
We study the equation $-\Delta u+u^q=0$, $q>1$, in a bounded $C^2$ domain $\Omega\subset \BBR^N$.
A positive solution of the equation is \emph{moderate} if it is dominated by a harmonic function and \emph{$\sigma$-moderate} if it is the limit of an increasing sequence of moderate solutions.  It is known that in the subcritical case, $1<q<q_c=(N+1)/(N-1)$, every positive solution is $\sigma$-moderate \cite{MVsub}. More recently Dynkin proved, by probabilistic methods,  that this remains valid in the supercritical case for $q\leq 2$, \cite{Dbook2}. The question remained open for $q>2$. In this paper we prove that, for all $q\geq q_c$, every positive solution is $\sigma$-moderate. We use purely analytic techniques which apply to the full supercritical range. The main tools come from linear and non-linear potential theory.
Combined with previous results, this establishes a 1-1 correspondence between positive solutions and their boundary traces in the sense of  \cite{MVfine}.
\end{abstract}
\maketitle
\section{Introduction}
In this paper we study \bvp{s} for the equation
\begin{equation}\label{q-eq}
    -\Gd u+|u|^q\sign u=0,\q q>1
\end{equation}
in a bounded $C^2$ domain $\Gw$. We say that $u$ is a solution of this equation if  $u\in L^q\loc(\Gw)$ and the equation holds in the sense of distributions. Every solution of the equation is in $W^{2,\infty}\loc(\Gw)$. In particular, every solution is in $C^1(\Gw)$.

Let $\GTM(\bdw)$ denote the space of finite Borel measures on the boundary. Put $$\gr(x):=\dist(x,\bdw)$$
and denote by $L^q_\gr(\Gw)$ the Lebesgue space with weight $\gr$.

For $\nu\in \GTM(\bdw)$ a (classical) weak solution of the \bvp
\begin{equation}\label{q-bvp}
    -\Gd u+|u|^q\sign u=0 \txin{\Gw}, \q u=\nu\txon{\bdw}
\end{equation}
is a function $u\in L^1(\Gw)\cap  L^q_\gr(\Gw)$ \sth
\begin{equation}\label{inteq}
   -\int_\Gw u\Gd \gf\,dx+\int_\Gw |u|^q\sign u\gf\,dx=-\int_{\bdw}\prtn\gf\,d\nu,
\end{equation}
for every $\gf\in C_0^2(\bgw)$ where
\begin{equation}\label{C02}
   C_0^2(\bgw):=\{\gf\in C^2(\bgw):\, \gf=0\txon{\bdw}\}.
\end{equation}

The \bvp \req{q-bvp} with data given by a finite Borel measure is well understood. It is known that, if a solution exists then it is unique. Gmira and V\'eron \cite{GV} proved that, if $1<q<(N+1)/(N-1)$, the problem possesses a solution for every $\nu\in \GTM(\bdw)$; if $q\geq (N+1)/(N-1)$ then the problem has no solution for any measure $\nu$ concentrated at a point.
The number $q_c:=(N+1)/(N-1)$ is  \emph{the critical value} for \req{q-bvp}. The interval $(1,(N+1)/(N-1))$ is  the \emph{subcritical} range ; the interval $[(N+1)/(N-1),\infty)$ is the \emph{supercritical} range.

In the early 90's the \bvp \req{q-bvp} became of great interest due to its relation to branching processes and superdiffusions (see Dynkin \cite{Dy91, Dy93}, Le Gall \cite{LG93a}). At first, the study of the problem concentrated on the characterization of the family of finite measures for which \req{q-bvp} possesses a solution. This question is closely related to the characterization of removable boundary sets. A compact set $K\sbs \bdw$ is \emph{removable} if every positive solution $u$ of \req{q-eq} which has a continuous extension to $\bgw\sms K$ can be extended to a function in $C(\bgw)$.


In a succession of works by Le Gall \cite{LG93b,LG95} (for $q=2$), Dynkin and Kuznetsov \cite{DK96,DK98a} (for $1<q\leq 2$) and Marcus and V\'eron \cite{MVsuper,MVrem} (the first for $q\geq 2$, the second providing a new proof for all $q\geq q_c$) the following results were established.
\bsn{Theorem A}  Let $K$ be a compact subset of $\bdw$. Then
\begin{equation}\label{rem-cr}
  K \text{  is removable }\iff C^{2/q,q'}(K)=0.
\end{equation}
 Here $q'=q/(q-1)$ and $C^{2/q,q'}$ denotes Bessel capacity on $\bdw$.
\esn

\bsn{Theorem B} Let $\nu\in \GTM(\bdw)$.  Problem \req{q-bvp} possesses a solution if and only if
$\nu\prec C^{2/q,q'}$, i.e. $\nu$ vanishes on every Borel set $E\sbs\bdw$ \sth $C^{2/q,q'}(E)=0$.
\esn
\note{Remark A.1} For solutions in  $L^q_\gr(\Gw)$, the removability criterion applies to signed solutions as well.
\vskip 2mm
\note{Remark A.2} For a non-negative solution $u$ of \req{q-eq}, the removability criterion can be extended to an arbitrary set $E\sbs \bdw$. Suppose that $u$ vanishes on every \qop subset of $\bdw\sms E$. Then
$$C^{2/q,q'}(E)=0\Lra u=0.$$
 This is a consequence of the capacitary estimates of \cite{MVcapest}.

 In view of the estimates of Keller \cite{Kell} and Osserman \cite{Oss} equation \req{q-eq} possesses solutions which are not in $L^q_\gr(\Gw)$. In particular the equation possesses solutions  which blow up everywhere on the boundary (recall that we assume that $\Gw$ is of class $C^2$).
 Such solutions, called \emph{large solutions} have been studied for a long time (see e.g. Loewner and Nirenberg \cite{LN} who studied the case $q=(N+2)/(N-2)$). It was established that the large solution is unique and its asymptotic behavior at the boundary was described
  (see  Bandle and Marcus \cite{BM90,BM95} and the references therein). The uniqueness of large solutions was also established for domains of class $C^0$ and even for \qop sets (see Marcus and V\'eron \cite{MVunique, MVcapint}).

  The next question in the study of equation \req{q-eq} was whether it is possible to assign to arbitrary solutions a measure, not necessarily finite, which uniquely determines the solution.
   (Eventually such a measure was  called a \emph{boundary trace.})
  In investigating this question, attention was restricted to positive solutions. The Herglotz theorem for positive harmonic functions served as a model. But, in contrast to the linear case, here one must allow unbounded measures.

  In \cite{LG93b} Le Gall studied \req{q-eq} with $q=2$ and $\Gw$ a disk in $\BBR^2$. He showed that, in this case, every positive solution possesses a boundary trace which uniquely determines the solution. The boundary trace was described in probabilistic terms and the proof relied mainly on probabilistic techniques.

 In \cite{MVcras96} Marcus and V\'eron introduced a notion of boundary trace (later  Dynkin called it
 `the rough trace') which can be described as a (possibly unbounded) Borel measure $\nu$
 with the following properties. There exists a closed set $F\sbs \bdw$ \sth

 (i) $\nu(E)=\infty$ for every non-empty Borel subset of $F$,

 (ii) $\nu$ is a Radon measure on $\bdw\sms F$.
\vskip 2mm

Let us denote the family of positive measures possessing these properties by $\CB_{\rm reg}(\bdw)$.
Given a positive solution $u$ of \req{q-eq}, we say that it has (rough)  boundary trace $\nu\in \CB_{\rm reg}(\bdw)$  if (with $F$ as above)
\vskip 1mm

(i') For every open \ngh $Q$ of $F$, $u\in L^1(\Gw\sms A)\cap L^q_\gr(\Gw\sms A)$ and
\req{inteq} holds for every $\vgf\in C^2_0(\bar \Gw)$  vanishing in a \ngh of $F$.
\vskip 1mm

(ii') If $\gx\in F$ then, for every open \ngh $A$ of $\gx$,
$$\int_{A\cap\Gw} u^q\gr\,dx=\infty.$$
\vskip 1mm

 The following result (announced in \cite{MVcras96}) was proved in \cite{MVsub}:

\bsn{Theorem C} Every positive solution of \req{q-eq} possesses a boundary trace in $\CB_{\rm reg}(\bdw)$.

If $1<q<q_c$ then, for every $\nu\in\CB_{\rm reg}(\bdw)$,  \req{q-eq} possesses a unique solution with boundary trace $\nu$.
\esn

In the supercritical case  it was  shown  in \cite{MVsuper} that, under some additional conditions on $\nu$, -- mainly that $\nu$ must vanish on subsets of $\bdw\sms F$ of $C^{2/q,q'}$-capacity zero, --
\req{q-eq} possesses a solution with rough trace $\nu$. These conditions were shown to be necessary and sufficient for existence. However, it soon became apparent that in the supercritical case,  the solution  is no longer unique. A counterexample to this effect was constructed by Le Gall in 1997. Therefore, in order to deal with the super critical case,  a more refined definition of boundary trace was necessary.

Kuznetsov \cite{Kz98} and Dynkin and Kuznetsov \cite{DK98b} provided such a definition, which they called `the fine trace'.
Their definition was similar to that of the rough trace, but the singular set $F$ was not required to be closed in the Euclidean topology. Instead it was required to be  closed  \wrto a finer topology defined in probabilistic terms. With this definition they showed that, if $q\leq 2$ then, for  any positive `fine trace' $\nu$,  \req{q-eq} possesses a solution the trace of which is equivalent, but not necessarily identical, to $\nu$. The equivalence is defined  in terms of polarity. Furthermore they showed that the minimal solution corresponding to a prescribed trace is $\gs$-moderate and it is \emph{the unique solution in this class.} The restriction to $q\leq 2$ is due to the fact that the proof was based on probabilistic techniques which do not apply to $q>2$.

A $\gs$-moderate solution was defined as the limit of an increasing sequence of positive moderate solutions. We recall that a moderate solution is a solution in $L^1(\Gw)\cap L^q_\gr(\Gw)$, i.e., a solution whose boundary trace is a finite measure.

In around the year 2002, Mselati proved in his Ph.D. thesis (under the supervision of Le Gall) that \emph{for $q=2$ every positive solution of \req{q-eq} is $\gs$-moderate.} This work appeared in \cite{Ms}. Mselati used a combination of analytic and probabilistic techniques such as the 'Brownian snake'  developed by Le Gall \cite{LGbook}. Following this, Dynkin \cite{Dbook2} extended Mselati's result proving:
 \vskip 2mm
\emph{If $q_c\leq q\leq 2$ then every positive solution of  \req{q-eq} is $\gs$-moderate.}
\vskip 2mm
Instead of the 'Brownian snake' technique, which can be applied only to the case $q=2$, Dynkin's proof used
new results on Markov processes that are applicable  to $q\leq2$.


At about the same time Marcus and V\'eron introduced a notion of boundary trace --  they called it `the precise trace' --  based on the classical notion of  $C^{2/q,q'}$-fine topology (see \cite{AH}).
A Borel measure $\nu$ on $\bdw$ belongs to this family of traces, to be denoted by $\CF^{2/q,q'}(\bdw)$, if there exists a $C^{2/q,q'}$-finely closed set $F\sbs \bdw$ \sth:
\vskip 2mm

(i) $\nu(E)=\infty$ for every non-empty Borel subset of $F$.

 (ii) Every point $x\in \bdw\sms F$ has a \qop \ngh $Q_x$ \sth $\nu(Q_x)<\infty$.

 (iii) If $E$ is a Borel set \sth $\nu(E)<\infty$ then $\nu$ vanishes on subsets of $E$ of $C^{2/q,q'}$-capacity zero.
\vskip 2mm

 In the subcritical case the $C^{2/q,q'}$-fine topology is identical to the Euclidean topology and consequently the precise trace coincides with the rough trace.

With this definition they proved \cite{MVfine}, by purely analytic methods:

\bsn{Theorem D} For every $q\geq q_c$:

 $(a)$ Every positive solution of \req{q-eq} possesses a boundary trace $\nu\in \CF^{2/q,q'}(\bdw)$.

 $(b)$ For every measure $\nu\in \CF^{2/q,q'}(\bdw)$, problem \req{q-bvp} possesses a $\gs$-moderate solution.

 $(c)$ The solution is unique in the class of $\gs$-moderate solutions.
 \esn

The question whether every positive solution of \req{q-eq} with $q>2$ is $\gs$-moderate remained open.
In the present paper we settle this question proving,

\bsn{Theorem 1} For every $q\geq q_c$, every  positive solution of \req{q-eq} is $\gs$-moderate.
\esn

 The proof employs only analytic techniques and \emph{applies to all} $q\geq q_c$. Of course the statement is also valid in the subcritical case, in which case it is an immediate consequence of Theorem {\bf C}.

 Combining Theorems {\bf C}, {\bf D} with Theorem {\bf 1} we obtain:
\bsn{Corollary 1}  For every $q>1$ and every non-negative $\nu \in   \CF^{2/q,q'}(\bdw)$,
problem \req{q-bvp} possesses a  unique solution. If $1<q<q_c$, $\CF^{2/q,q'}(\bdw)=\CB_{\rm reg}(\bdw)$.
\esn


The method developed in the present paper can be adapted and applied  to  a general class of problems which includes \bvp{s} for equations such as
$$-\Gd u+\gr^\ga |u|^q\sign u=0,\q \ga>-2$$
and
$$-\Gd u+g(u)=0,$$
where $g\in C(\BBR)$ is odd, monotone increasing and satisfies the $\Gd_2$ condition and the Keller--Osserman condition. For equations of the latter type, the method can be adapted to \bvp{s} in  Lipschitz domains as well. These results will be presented in a subsequent paper.

\vskip 2mm

The main ingredients used in the present paper are:
\vskip 2mm

 (a) Nonlinear potential theory and fine topologies associated with Bessel capacities (see  \cite{AH} and \cite{MVfine}).
 \vskip 1mm

 (b) The theory of  boundary value problems for equations of the form
 $$L^Vu:=-\Gd u+Vu=0 \q\text{in}\;\Gw,$$
 where $V>0$ and $\gr^2 V$ is bounded. Here we use mainly the results of Ancona \cite{An87} together with classical potential theory results (see e.g.  \cite{An-SLN}).

 \vskip 3mm

 \noindent\textbf{Acknowledgment} The author wishes to thank Professor Alano Ancona for several most helpful discussions. He also wishes to thank him for two very useful personal communications \cite{An-05} (later published in \cite{An-App}) and \cite{An-10}.

\section{Preliminaries: on the equation $-\Gd u+Vu=0$.}

For the convenience of the reader we collect here some definitions and results of classical potential theory concerning  operators of the form $L^V=-\Gd +V$, that will be used in the sequel. The results apply also to operators of the form $-L_0+V$ where $L_0$ is a second order uniformly elliptic operator on  differentiable manifolds with negative curvature. However we shall confine ourselves to the operator $L^V$ in a bounded domain $\Gw\sbs \RN$ which is either a $C^2$ domain or Lipschitz.

The following conditions on $V$ will be assumed, without further mention, throughout the paper.
\begin{equation}\label{V-profile}
  0<V\leq c\gr (x)^2,\q V\in C(\Gw).
\end{equation}

By \cite{An87}  if $\Gw$ is a bounded Lipschitz domain, the Martin boundary can be identified with the Euclidean boundary $\bdw$ and, for every $\gz\in \bdw$ there exists a positive $L^V$ harmonic which vanishes on $\bdw\sms\{\gz\}$. If normalized this harmonic is unique. We choose a fixed reference point, say $x_0\in \Gw$ and denote by
$K_\gz^V$ this $L^V$ harmonic, normalized by $K_\gz^V(x_0)=1$.

We observe that the positivity of $V$ is essential for this result. Indeed the result depends on the weak coercivity of $L^V$ (see definition in \cite[Section 2]{An87}) which is guaranteed in our case by Hardy's inequality.

As a consequence of the above one obtains the following basic result (see Ancona \cite{An87}, Theorem 3 and  Corollary 13),

\bsn{Representation Theorem}
For each positive $L^V$-harmonic function $u$ in $\Gw$ there exists a unique positive measure $\mu$ on  $\bdw$ \sth
\begin{equation}\label{martin}
  u(x)=\int_{\bdw} K^V_\gz d\mu(\gz) \forevery x\in \Gw.
\end{equation}
\esn

The function
$$K^V(\cdot,\gz)=K^V_\gz(\cdot)$$
is the Martin kernel. In $C^2$-domains, \wrto a classical elliptic operator such as $-\Gd$, it can be identified with the Poisson kernel $P$. More precisely in this case
$$K^0(\cdot,\gz)=P(\cdot,\gz)/P(x_0,\gz),$$
where $x_0$ is a fixed reference point in $\Gw$. In Lip domains, \wrto $-\Gd$, $K^0$  is precisely the harmonic measure.

In the sequel  we write $$K_\gz:=K^0_\gz.$$ The measure $\mu$ corresponding to an $L^V$ harmonic $u$ will be called \emph{the $L^V$ boundary trace of $u$} and we use the notation
\begin{equation}\label{martin'}
  \BBK^V_\mu:=\int_{\bdw} K^V_\gz d\mu(\gz), \q \BBK_\mu:=\int_{\bdw} K_\gz d\mu(\gz).
\end{equation}

Let $D$ be a Lipschitz domain \sth $\bar D\sbs\Gw$ and $h\in L^1(\prt D)$. We denote by $S^V(D, h)$
the solution of the problem
\begin{equation}\label{LV-bvp}
 L^Vw:=-\Gd w+Vw=0 \text{ in $D$},\q w=h \text{ on $\prt D$.}
\end{equation}
If $\mu$ is a finite measure on $\prt D$, $S^V(D, \mu)$ is defined in the same way. If $D$ is a $C^2$ domain, a function $w\in L^1(D)$ is a solution of \req{LV-bvp} (with $h$ replaced by $\mu$) if
\begin{equation}\label{LV-bvp-mu}
\int_D(-w\Gd\vgf+Vw\vgf)dx=-\int_{\prt D} \prtn\vgf \,d\mu,
\end{equation}
for every $\vgf\in C^2_0(\bar D)$.

A family of domains $\{\Gw_n\}$ \sth $\bgw_n\sbs \Gw_{n+1}$ and $\cup \Gw_n=\Gw$ is called an \emph{exhaustion of $\Gw$. } We say that $\{\Gw_n\}$ is a Lipschitz (resp. $C^2$) exhaustion if each domain $\Gw_n$ is Lipschitz (resp. $C^2$).

An l.s.c. function $u\in L^1\loc(\Gw)$ is $L^V$ \superh if $L^Vu\geq 0$ in distribution sense. Such a function is necessarily in $W^{1,p}\loc(\Gw)$ for some $p>1$ and \consy it possesses an $L^1$ trace on $\prt D$ for every $C^2$ domain $D\Subset \Gw$. Furthermore, for every such domain, $u\geq S^V(D,u)$. If $u$ is positive, the same holds for every Lipschitz domain $D\Subset\Gw$.

If $u$ is an $L^V$-\superh in $\Gw$ and $D$ a $C^2$ domain \sth $D\Subset \Gw$ then
the function $u\ind{D}$ defined by
$$u\ind{D}=S^V(D,u) \text{ in $D$,}\q u\ind{D}=u \text{ in $\Gw\sms D$}$$
is called the \emph{$D$-truncation of $u$.} This function is an $L^V$-\superh.

\blemma{M0} Let $u$ be a non-negative $L^V$-\superh and $\{\Gw_n\}$ a Lipschitz exhaustion of $\Gw$. Then the following limit exists
\begin{equation}\label{sup-minorant}
   \tl u:=\lim S^V(\Gw_n, u)
\end{equation}
and $\tl u$ is the largest $L^V$-harmonic dominated by $u$.
\es
\bproof The \seq $\{S^V(\Gw_n, u)\}$ is non-increasing. \Consy the limit exists and it is an $ L^V$ -harmonic. Every $L^V$ harmonic $v$ dominated by $u$ must satisfy $v\leq S^V(\Gw_n, u)$ in $\Gw_n$. Therefore $\tl u$ is the largest such harmonic.
\eproof

\bdef{Potential} A non-negative $L^V$-superharmonic  is called an $L^V$- potential if its largest $L^V$-harmonic minorant is zero.
\es

The following is an immediate consequence of \rlemma{M0}.
\blemma{M1} A non-negative superharmonic function $p$ is an $L^V$ potential \ifif
$$S^V(D_\gb,p)\to 0 \txt{as} \gb\to0.$$
\es

\bsn{Riesz decomposition theorem} Every non-negative $L^V$-superharmonic $u$ can be written in a unique way in the form $u=p+h$ where $p$ is an $L^V$ potential and $h$ a non-negative $L^V$-harmonic.
\esn

\note{Remark}In fact $h=\tl u$ as defined in \req{sup-minorant}.

For further results concerning the $L^V$-potential see \cite[Ch.I, sec. 4]{An-SLN}.

\bdef{R^A_s} Let $A\sbs \Gw$ and let $s$ be a positive $L^V$-superharmonic.  Then $R^A_s$
(called the reduction of $s$ relative to $A$) is given by
$$R_s^A=\textnormal{lower envelope of }\{f:\, 0\leq f \text{ superharmonic,}\; s\leq f \text{ on $A$.}\}$$
\es

If $A$ is open then $R_s^A$ itself is $L^V$-superharmonic so that the lower envelope is simply the minimum, \cite[p.13]{An-SLN}.

\bdef{thin_set} Let $\gz\in \bdw$.
A set $A$ is $L^V$ thin at $\gz$ (in French '$A$ est $\gz$-effil\'{e}') if $R_{K^V_\gz}^A\not \equiv K^V_\gz$.
\es
In view of a theorem of Brelot, if $A$ is open:
$$R_{K^V_\gz}^A\not \equiv K^V_\gz \iff R_{K^V_\gz}^A \text{ is an $L^V$-potential.}$$

Furthermore, even if $A$ is not open there exists an open set $O$ \sth $A\sbs O$ and $O$ is thin at $\gz$.

\blemma{M2} Assume that $A$ is thin at $\gz\in \bdw$ and $A$ open. Let $\{D_n\}$ be a $C^2$ exhaustion of $\Gw$ and put $A_n=\bdw_n\cap A$. Then
$$S^V(\Gw_n, K^V_\gz\chi\ind{A_n})\to0.$$
\es

\bdef{Vfine} Let $\gz\in \bdw$ and $f$ a real function on $\Gw$. We say that $f$ admits the fine limit $\ell$ at $\gz$ if there exists a closed set $E\sbs \Gw$ \sth $E$ is thin at $\gz$ and
$$\lim_{x\to\gz,\,x\in \Gw\sms E}f(x)=\ell.$$
To indicate this type of convergence we write,
$$\lim_{x\to\gz}f(x)=\ell,\; L^V-\text{finely}.$$
\es

Recall that there also exists an open set $A$ \sth $E\sbs A$ and $A$ is thin at $\gz$.

\bprop{fine-nt}
Let $u$ be a positive  $L^V$ harmonic function, or  a solution of \req{q-eq}. For $\gz\in \bdw$,
$$\lim_{x\to\gz} u=b\q\text{$L^V$- finely} \Lra \lim_{x\to\gz} u=b \text{ n.t.,}$$
where  'n.t.' means 'non-tangentially'.
\es
\bproof Let $\gr(x):=\dist(x,\bdw)$. By \cite[Lemma 6.4]{An-SLN}, if $A$ is an $L^V$ thin set at $\gz$ and $\gb_n\dar 0$ then
$$A\cap\{x\in \Gw:\, |x-\gz|<\gr(x),\; \gb_n/2<\gr(x)<3/2\gb_n\}\neq\ems$$
for all sufficiently large $n$. Therefore  the assertion follows from Harnack's inequality.
\eproof

For the next two theorems see \cite[Prop.1.6 \& Thm. 1.8] {An-SLN}.
\bth{Fatou0} If $p$ is an $L^V$ potential then, for every positive $L^V$ harmonic $v$:
$$\lim_{x\to\gz,\,fine} p/v=0 \q\mu_v - a.e.$$
where $\mu_v$ is the $L^V$ boundary trace of $v$.
\es

\bth{Fatou} \textnormal{\textbf{[Fatou-Doob-Naim]}} If $u,v$ are positive $L^V$ harmonics then $u/v$ admits a fine limit $\mu_v$ a.e. Furthermore
$$\lim_{x\to\gz,\,fine} u/v= f=\frac{d\mu_u}{d\mu_v} \q \mu_v-a.e.$$
where $\mu_u$ and $\mu_v$ are the $L^V$ boundary traces of $u$ and $v$ respectively and the term on the \rhs denotes the Radon-Nikodym derivative.
\es

The next lemma -- an application of the theorem of Fatou -- is due to Ancona \cite{An-10}.
\blemma{v_ con-nt} Assume that $v$ is a positive $L^V$ harmonic function with $L^V$ boundary trace $\nu$.
 Then
\begin{equation}\label{lim_v>0}
   \lim_{x\to\gz} v>0 \q\text{n.t. $\nu$-a.e.} \gz\in \bdw.
\end{equation}
If $\nu\perp\BBH_{N-1}$ then
\begin{equation}\label{lim_v=in}
   \lim_{x\to\gz} v=\infty \q\text{n.t. $\nu$-a.e. }
\end{equation}
\es
\bproof   The function $1$ is an $L^V$ \superh. If it is a potential then, by \rth{Fatou0}
$$\lim_{x\to\gz} 1/v=0\; \text{ $L^V$-finely $\nu$-a.e.}$$
Therefore, by \rprop{fine-nt}
$$\lim_{x\to\gz}v=\infty \q\text{n.t. $\nu$-a.e.}$$
If $1$ is not a potential there exists a positive $L^V$ harmonic $w$ and a potential $p$ \sth $1=w+p$ . Let $w=K^V_\gg$ and put $d\gg/d\nu=:f$.
 By \rth{Fatou}
$$\lim w/v=f \; \text{ $L^V$-finely $\nu$-a.e.}$$
(We do not exclude the possibility that $f=0\; \nu$-a.e. but, of course, $f<\infty\;\nu$-a.e. )
Since $p/v\to 0$ finely $\nu$-a.e., it follows that
$$\lim 1/v=\lim (w+p)/v=f\; \text{ $L^V$-finely $\nu$-a.e.}$$
Applying again \rprop{fine-nt} we obtain
$$\lim 1/v=f \text{ n.t. $\nu$-a.e.}$$
which in turn implies \req{lim_v>0}.

If $\nu\perp \BBH_{N-1}$ then $f=0$ $\nu$-a.e. and \consy $v\tin$ n.t. $\nu$-a.e.
\eproof


\section{Moderate  solutions of $L^Vu=0$}
We recall some definitions from  \cite{Dbook2} following the notation of \cite{An-App}.
\bdef{LV-regular}
We shall say that a boundary point $\gz$ is $L^V$ regular if $\tl K^V(\cdot,\gz)$ $(=$ the largest $L^V$ harmonic dominated by the $L^V$ \superh $K(\cdot,\gz))$ is positive. The point $\gz$ is $L^V$ singular if
$\tl K^V(\cdot,\gz)=0$, i.e. $K(\cdot,\gz)$ is an $L^V$ potential.

We denote by $Sing(V)$, respectively $Reg(V)$, the set of singular points, respectively regular points, of $L^V$.
\es

\Remark
If $\gz$ is $L^V$-regular then
$$\tl K^V(\cdot,\gz)=c(x_0) K^V(\cdot,\gz), \q c(x_0)=\tl K^V(x_0,\gz), $$
where $x_0$ is a fixed reference point in $\Gw$ \sth
$$K^V(x_0,\gz)=1 \forevery \gz\in \bdw.$$

\note{Notation} The family of finite Borel measures on a set $A$  is denoted by $\GTM(A)$. For $A=\bdw$ we shall write simply $\GTM$. If $\mu\in \GTM(A)$ we denote by $|\mu|$ the total variation measure and by $\norm{\mu}\ind{\GTM(A)}$ the total variation norm.


\bdef{Vmod} 

An $L^V$ harmonic $u$ is $L^V$-moderate if $u\in L^1(\Gw)\cap L^1(\Gw;V\gr)$ and there exists a measure $\nu\in \GTM$ \sth
\begin{equation}\label{eqint}
     \int_\Gw(-u\Gd\vgf+ uV\vgf)dx=-\int_{\bdw}\prtn\vgf\,d\nu,
\end{equation}
for every $\vgf\in C_0^2(\bar\Gw)$.

A  measure $\nu\in \GTM$ is $L^V$-moderate if there exists a moderate $L^V$ harmonic satisfying \req{eqint}.

The space of $L^V$ moderate measures is denoted by $\GTM^V$.
\es

\bdef{m-trace} Let $u\in W^{1,p}\loc(\Gw)$ for some $p>1$. We say that $u$ possesses an m-boundary trace $\nu\in \GTM(\bdw)$ if, for every $C^2$-exhaustion of $\Gw$, say $\{\Gw_n\}$,
$$u\lfloor\ind{\bdw_n}d\BBH_{N_1}\rightharpoonup \nu,$$
weakly \wrto $C(\bgw)$, i.e.,
\begin{equation}\label{m-trace}
  \int_{\bdw_n}uh\,dS\to \int_{\bdw}h\,d\nu  \forevery h\in C(\bar\Gw).
\end{equation}
If $\nu$ is the m-boundary trace of $u$ we write $\tr\,u:=\nu$.
\es
\vskip 1mm

\note{Remark} If $u$ possesses an m-boundary trace $\nu$ then $u\in L^1(\Gw)$ and
\begin{equation}\label{|u|dS}
 \sup \int_{\bdw_n}|u| dS<\infty.
\end{equation}
This follows immediately from the definition. It is easily verified that, if $u$ is $L^V$ moderate and satisfies \req{m-trace} then $\nu$ is the m-boundary trace of $u$.

\vskip 2mm

\note{Notation} Let $\gr(x):=\dist(x,\bdw)$. In the case of $C^2$ domains, there exists $\gb_0>0$ \sth for $x\in \Gw$, $\gr(x)<\gb_0$, there exists a unique point on $\bdw$, to be denoted by  $\gs(x)$, \sth
$$|x-\gs(x)|=\gr(x).$$
Thus
$$ x-\gs(x)=\gr(x)\gn\ind{\gs(x)}$$
where $\gn_\gz$ denotes the unit normal at $\gz\in \bdw$ pointing into the domain.
(We also denote $-\gn_\gz=:\mathbf{n}_\gz$.) It can be shown that the function $x\mapsto \gs(x)$ is in $C^2(\bar\Gw^0)$ where
$$\Gw^0:=\{x\in \Gw:\,0<\gr(x)<\gb_0\}.$$
The mapping $x\mapsto (\gr(x),\gs(x))$ is a $C^2$ homeomorphism of $\Gw^0$ onto
 $$\{(\gr,\gs)\in \BBR_+\ti \bdw: 0<\gr<\gb_0\}.$$
Thus $(\gr,\gs)$  can be used as an alternative set of coordinates in $\Gw^0$; we call them `flow coordinates'.

Put,
$$D_\gb=[x\in \Gw,\;\gr(x)>\gb],\q \Gw_\gb=\Gw-\bar D_\gb,\q \Gs_\gb=[x\in\Gw,\;\gr(x)=\gb].$$
and for $\ga\in (0,\infty)$ and  $\gz\in \bdw$
$$\CC^\ga_\gz= \{x\in \Gw^0: |x-\gz|>\ga\gr(x)\}.$$
When $\ga=1$,  the upper index will be omitted.

In the sequel we assume that \emph{$\Gw$ is a bounded $C^2$ domain.}

\blemma{mod-basic}

$(i)$ If $v$ is a positive $L^V$  \superh and
\begin{equation}\label{LVr}
    \int_\Gw vV\gr\,dx<\infty.
\end{equation}
then $v$ is moderate. In particular, $v\in L^1(\Gw)$ and it possesses an m-boundary trace
$\nu\in \GTM$. The supremum of $L^V$-harmonics dominated by $v$, say $v'$, is an $L^V$ harmonic and has the same m-boundary trace.

$(ii)$ If $v$ is a positive $L^V$ superharmonic  and $v$ possesses an  m-boundary trace
$\nu\in \GTM$ then $v'$, the supremum of $L^V$-harmonics dominated by $v$, is an $L^V$ moderate harmonic and $\tr\,v'\le \nu$. If $v$ is not a potential then $v'$ is positive.

$(iii)$ If $v$ is an $L^V$ harmonic (not necessarily positive) and $v$ possesses an m-boundary trace $\nu$ then $v$ is $L^V$ moderate.
\es

\bproof
(i) Let $w\in L^1(\Gw)$ be the (unique) solution of the problem
\begin{equation}\label{Gdw}
   -\Gd w=Vv \text{ in $\Gw$}, \q w=0 \text{ on $\bdw$}.
\end{equation}
The solution exists because $v$ satisfies \req{LVr}. Then
$$-\Gd(w+v)\geq 0$$
and \consy $w+v\in L^1(\Gw)$ and it possesses an m-boundary trace $\nu\in \GTM(\bdw)$. As $w\in L^1(\Gw)$ and has  m-boundary trace zero, it follows that $v\in L^1(\Gw)$ and has m-boundary trace $\nu$.

Given $\vgf\in C_0^2(\bar\Gw)$, for each $\gb\in (0,\gb_0/2)$ we can construct a function $\vgf_\gb\in C^2_0(\bar D_\gb)$ \sth
$$\vgf_\gb(x)=\vgf(\gr(x)-\gb, \gs(x)) \text{ for } \gb\leq \gr(x)<\gb+\gb_0/4$$
and
$$\vgf_\gb\to\vgf \text{ in }C^2(\Gw),\q \sup_\gb\norm{\vgf_\gb}\ind{C^2(\bar D_\gb)}<\infty.$$
Let
$$v_\gb=S^V(D_\gb,v).$$
Then $0\leq v_\gb\leq v$ and $v_\gb\dar v'$ as $\gb\dar 0$. Furthermore
$$     \int_{D_\gb}(-v_\gb\Gd\vgf_\gb+ v_\gb\,V\vgf_\gb)dx=-\int_{\prt D_\gb}\prtn\vgf_\gb\,v\,dS.$$
Since $v\in L^1(\Gw;V\gr)\cap L^1(\Gw)$ and $v$ possesses m- boundary trace $\nu$ on $\bdw$ we obtain (by going to the limit as $\gb\to0$):
\begin{equation}\label{eqint'}
     \int_\Gw(-v'\Gd\vgf+ v'V\vgf)dx=-\int_{\bdw}\prtn\vgf\,d\nu,
\end{equation}
for every $\vgf\in C_0^2(\bar\Gw)$.

\medskip

(ii) If $v$ is a positive $L^V$ superharmonic then there exists a Radon measure $\gl>0$ in $\Gw$ \sth $L^Vu=\gl$ in the sense of distributions. Therefore $v\in W^{1,p}\loc(\Gw$ for some $p>1$. Let $v_\gb=S^V(D_\gb, v)$. Then $v_\gb\le v$ in $D_\gb$ and
\begin{equation}\label{vgb1}
     \int_{D_\gb}(-v_\gb\Gd\vgf+ v_\gb\,V\vgf)dx=-\int_{\prt D_\gb}\prtn\vgf\,v\,dS
\end{equation}
for every $\vgf\in C^2(\bar D_\gb)$. Choosing $\vgf$ to be the solution of
$$-\Gd \vgf=1 \text{ in $D_\gb$},\q  \vgf=0 \text{ on $\prt D_\gb$}$$
we obtain
\begin{equation}\label{vgb2}
\norm{v_\gb}\ind{L^1(D_\gb)}+\norm{v_\gb}\ind{L^1(D_\gb;V\gr_\gb)}\leq C \int_{\prt D_\gb}v\,dS
\end{equation}
with a constant $C$ independent of $\gb$. Here $\gr_\gb$ is the first eigenfunction of $-\Gd$ in $D_\gb$ normalized by $\gr_\gb(x_0)=1$. Since, by assumption, $v$ has an m-boundary trace, the \rhs of the inequality is bounded. In addition, $\gr_\gb$ tends to the first normalized eigenfunction of $-\Gd$ in $\Gw$. Therefore $v_\gb\dar  v'$ as $\gb\dar 0$, locally uniformly in $\Gw$ and  $v'\in L^1(\Gw)\cap L^1(\Gw;V\gr)$. By \req{vgb1} with $\vgf=\vgf_\gb$ and Fatou's lemma we obtain -- using the fact that $v_\gb\leq v\in L^1(\Gw)$ and $\vgf_\gb\to\vgf$ in $C^2(\Gw)$ --
\begin{equation}\label{eqint<}
     \int_\Gw(-v'\Gd\vgf+ v'V\vgf)dx\leq-\int_{\bdw}\prtn\vgf\,d\nu,
\end{equation}
for every non-negative $\vgf\in C_0^2(\bar\Gw)$. \Consy (by a standard argument)  $v'$ has an m-boundary trace, say $\nu'$,
\sth $\nu'\le \nu$.

If, in addition, $v$ is not an $L^V$ potential then $v'>0$.

(iii) The proof is essentially the same as that of part (ii) except that, in the present case,
 inequality \req{vgb2} is replaced by
\begin{equation}\label{vgb2'}
\norm{v_\gb}\ind{L^1(D_\gb}+\norm{v_\gb}\ind{L^1(D_\gb;V\gr_\gb)}\leq C \int_{\prt D_\gb}|v|\,dS.
\end{equation}
This inequality is proved by a standard argument as in e.g. \cite{MVsuper}. Since $v$ is an $L^V$ harmonic, $v_\gb=v$ in $D_\gb$. Therefore we obtain $v\in L^1(\Gw)\cap L^1(\Gw;V\gr)$ and
$$     \int_{D_\gb}(-v\Gd\vgf_\gb+ vV\vgf_\gb)dx=-\int_{\prt D_\gb}\prtn\vgf_\gb\,v\,dS.$$
Finally, taking the limit as $\gb\to 0$, we obtain
$$   \int_{\Gw}(-v\Gd\vgf+ vV\vgf)dx=-\int_{\prt D_\gb}\prtn\vgf_\gb\,d\nu.$$
\eproof

\blemma{modsol} $(i)$ If $\nu\in \GTM^V$ then the solution of \req{eqint} is unique. It will be denoted by $\BBM^V_\nu$.

\noindent$(ii)$ The space $\GTM^V$ is linear and
\begin{equation}\label{MV-monotone}
0\leq \nu\iff 0\leq\BBM^V_\nu
\end{equation}

\noindent$(iii)$ Let $\tau\in\GTM$ and $\nu\in \GTM^V$.
\begin{align}\label{property1}
  |\tau|\leq \nu &\Lra\tau\in \GTM^V,\\
  \nu\in \GTM^V&\Lra |\nu|\in \GTM^V. \label{property2}
\end{align}

\noindent$(iv)$ If $\nu\in \GTM^V$ then
$|\BBM^V_\nu|\leq\BBM^V_{|\nu|}.$
\es
\bproof  (i), (ii) and \req{property1} are  classical. We turn to the proof of \req{property2}.
Put $u=\BBM^V_\nu$. Since $uV\in L^1_\gr(\Gw)$,  there exist $v^+$ and $v^-$ in $L^1(\Gw)\cap L^1(\Gw;V\gr)$ \sth
$$-\Gd v^\pm+Vu_{\pm}=0 \text{ in }\Gw,\q v^\pm=\nu_{\pm} \text{ on }\bdw.$$
It follows that
$$u=v^+ - v^-,\q |u|\leq v^+ + v^-=:w$$
and
$$-\Gd w +Vw\geq -\Gd w +V|u|=0.$$
Thus $w$ is a positive $L^V$ \superh with m-boundary trace $|\nu|$ and
$$w\in L^1(\Gw)\cap L^1(\Gw;V\gr).$$
By \rlemma{mod-basic} (i), the largest $L^V$ harmonic dominated by $w$, say $w'$, is $L^V$ moderate and
$$\tr\,w'=\tr\,w=|\nu|.$$
Thus $w'=\BBM_{|\nu|}$. Since $0<w'$ and $ u<w'$ it follows that $|u|\leq w'$, which is precisely assertion (iv).
\eproof

Denote,
\begin{equation}\label{GTMV0}
  \GTM^V_0:=\{\nu\in \GTM:\, \BBK_{|\nu|}V\in L^1_\gr(\Gw)\}.
\end{equation}

The following is an immediate consequence of \rlemma{mod-basic}:
\blemma{V0subV} $\GTM^V_0\sbs \GTM^V$.
\es

\bproof If $\nu\in\GTM^V_0$ then $\BBK_{|\nu|}$ is an $L^V$  \superh satisfying \req{LVr}.
\eproof

\Remark  In general there may exist positive measures $\nu$ \sth $\BBK^V_\nu\in L^1(\Gw;V\gr)$  but $\BBK_\nu\not\in L^1(\Gw;V\gr)$.

\blemma{VinLp} If $V\in L^{q'}_\gr(\Gw)$ for some $q'>1$ then every positive measure in $W^{-2/q,q}$ belongs to $\GTM^V_0$.
\es
\bproof If $\nu\in W^{-2/q,q}$ then $\BBK_\nu\in L^q_\gr(\Gw)$ (see \cite[1.14.4.]{Triebel} or \cite{MVbesov}). Therefore
$V\BBK_\nu \in L^1_\gr(\Gw)$. If, in addition, $\nu\geq 0$ then $\nu\in \GTM^V_0$.
\eproof
\Remark  There are  signed measures $\nu\in W^{-2/q,q}$ \sth
$|\nu|\not \in W^{-2/q,q}$. Therefore in general $W^{-2/q,q}$ may not be contained in  $\GTM^V_0$.
\vskip 2mm

\bprop{nu_nu'}
Let $v$ be a positive, $L^V$ moderate harmonic with m-boundary trace $\nu$. Let $\nu'$ be the $L^V$ boundary trace of $v$, i.e.,
\begin{equation}\label{MV_KV}
  \BBM_\nu^V=\BBK_{\nu'}^V.
\end{equation}
Then
 \begin{equation}\label{nu_nu'}
\nu'(E)=0\iff \nu(E)=0.
 \end{equation}
Furthermore,
 \begin{equation}\label{KV_L1}
 K^V(\cdot,\gz)V\in L^1_\gr(\Gw) \q \nu'-a.e.
\end{equation}

Let $F\sbs \bdw$ be compact and denote
$$v_F:=\inf(v,\BBK_{\nu_F}), \q \nu_F:=\nu\chi_F.$$
Then $v_F$ is a moderate supersolution of $L^V$ and the largest $L^V$ harmonic dominated by $v_F$ is given by
\begin{equation}\label{(vF)'}
  (v_F)'=\BBK^{V}_{\nu'\chi_F}.
\end{equation}
\es
\bproof

Since $vV\in L^1_\gr(\Gw)$ it follows (by Fubini) that
$$\int_{\bdw}(\int_\Gw K^V(\cdot,\gz)V\gr \big)d\nu'(\gz)<\infty$$
which implies \req{KV_L1}.

Let $F$ be a compact subset of $\bdw$. If $\nu'(F)>0$ then $K^V_{\nu'_F}$  (as usual $\nu'_F:=\nu'\chi_F$) is a positive $L^V$ moderate harmonic which vanishes on $\bdw\sms F$ and is dominated by $v$. Therefore
$$K^V_{\nu'_F}\leq v_F,$$
 and \consy $\nu(F)>0$.

 Next we show that $\nu(F)>0$ implies $ \nu'(F)>0$. Since $\nu(F)>0$, $v_F$ is a positive $L^V$ \superh with m- boundary trace $\nu_F$. In addition $v_F\leq v\in L^1(\Gw;V\gr)$.  Therefore by \rlemma{mod-basic} (i),  if $(v_F)'$ is the largest $L^V$ harmonic dominated by $v_F$
 then $(v_F)'$ is $L^V$ moderate with m- boundary trace $\nu_F$.
 Thus  $$0<(v_F)'\leq v_F .$$

 On the other hand, the largest $L^V$ harmonic dominated by $v$ and vanishing on $\bdw\sms F$ is $\BBK^V\ind{\nu'\chi_F}$. It follows that
 $$0<(v_F)'= \BBK^V\ind{\nu'\chi_F}$$
which implies that $\nu'(F)>0$.
\eproof
\bprop{LVregular}  $(i)$ For every $\gz\in \bdw$,
\begin{equation}\label{KVmoderate}
 K^V(\cdot,\gz)V\in L^1_\gr(\Gw)  \iff \gz\text{ is $L^V$-regular}.
\end{equation}

\noindent$(ii)$ If
$\nu$ is a positive measure in $\GTM^V$ then
$K^V(\cdot,\gz)$ is $L^V$-moderate $\nu$-a.e. and
\begin{equation}\label{nu(Sing)}
  \nu(Sing(V))=0.
\end{equation}


\noindent$(iii)$ If $V\in L^{q'}(\gw)$ for some $q'>1$ then $C_{2/q,q'}$-a.e. point $\gz\in \bdw$ is $L^V$ regular. (Here $\rec{q}+\rec{q'}=1$.)
\es
\bproof
(i) Assume that  $K^V(\cdot,\gz)V\in L^1_\gr(\Gw)$.  Then, by \rlemma{mod-basic}, $K^V_\gz$ is $L^V$ moderate. Its m-boundary trace $\tau_\gz\in \GTM$ is concentrated at $\gz$.
Thus  $\tau_\gz=a(\gz)\gd_\gz$ for some $a>0$. It follows that $K^V_\gz$ is a subsolution of the \bvp
$$-\Gd z=0 \text{ in $\Gw$}, \q z=\tau_\gz\text{ on $\bdw$.}$$
Therefore
$$K^V(\cdot,\gz)\leq a(\gz)K(\cdot,\gz).$$
This implies that $\tl K_\gz^V>0$, i.e., $\gz$ is regular.

 Assume that $\gz$ is $L^V$ regular. Then, by definition, $K(\cdot,\gz)$ is not a potential and has  m-boundary trace $=\gd_\gz$. By \rlemma{mod-basic}(ii)
 the largest $L^V$ harmonic dominated by $K(\cdot,\gz)$, which we denote by
 $\tl  K^V(\cdot,\gz)$, is $L^V$ moderate and its boundary trace, say $\gl$, is a positive measure bounded by $\gd_\gz$ . By uniqueness of the positive, normalized $L^V$ harmonic vanishing on $\bdw\sms\{\gz\}$,
 $$K^V(\cdot,\gz)=\tl K^V(\cdot,\gz)/\tl K^V(x_0,\gz).$$
 Thus $K^V(\cdot,\gz)V\in L^1_\gr(\Gw)$.
 \medskip

 (ii) By \req{KV_L1} and \req{nu_nu'},
\begin{equation}\label{KVmod}
K^V(\cdot,\gz)V\in L^1_\gr(\Gw) \q \nu-a.e.
\end{equation}
 By (i), this implies the second assertion.

\medskip

(iii) In this case,  every positive measure $\nu\in W^{-2/q,q}(\bdw)$  is in $\GTM^V_0$ which is contained in $\GTM^V$. It follows that the set of singular points of $L^V$ must have $\Capq$ - capacity zero.
\eproof

\section{Preliminaries: on the equation $-\Gd u +u^q=0$}

In this section we collect  some definitions and known results on positive solutions of \req{q-eq}, that will be needed for the proof of the main result.

A basic concept in this theory is that of $\Capq$ - fine topology. For the general theory of  $C_{m,p}$ capacity and $C_{m,p}$-fine topology we refer the reader to \cite{AH}. For more special results, useful in our study, we refer the reader  to the  summary in \cite[Section 2]{MVfine}.

The closure of a set $A\sbs\bdw$ in $\Capq$-fine topology will be denoted by $\tl A$. We shall say that two sets $A, B$ are $\Capq$ equivalent (or briefly q-equivalent) if $\Capq(A\,\Gd\,B)=0$.

There exists a constant $c$ \sth for every set $A$
$$\Capq(\tl A)\leq c\Capq(A).$$

We  recall the definition of regular and singular boundary points of a positive solution $u$ of \req{q-eq}.
 A point $\gz\in \bdw$ is a \emph {$q$-regular point} of $u$ if there exists a $\Capq$ \ngh of $\gz$, say $O_\gz$ \sth
$$\int_{O_\gz\cap\Gw}u^q\gr\,dx<\infty.$$
$\gz$ is q-singular if it is not q-regular. The set of q-regular points is denoted by $\CR(u)$ and the set of singular points by $\CS(u)$. Evidently $\CR(u)$ is $\Capq$ open.

If $F$ is a \qcl subset of $\bdw$ then there exists an increasing \seq of compact subsets $\{F_n\}$ \sth $\Capq(F\sms F_n)\to 0$.

If $u$ is a positive solution of \req{q-eq} we say that it vanishes on a \qop set $O=\bdw\sms F$ if it is the limit  of an increasing \seq of positive solutions $\{u_n\}$ \sth $u_n\in  C(\bar\Gw\sms F_n)$
and $u_n=0$ on $\bdw\sms F_n$. The q-support of the boundary trace of $u$ -- denoted by $\suppq u$ -- is the complement of the largest \qop subset of $\bdw$ where $u$  vanishes.

Let $\nu\in \GTM$. We say that $u$ is a solution of the problem
\begin{equation}\label{bvp_nu}
  -\Gd u+|u|^q\sign u=0 \text{ in $D$},  u=\nu\text{ on $\prt D$}
\end{equation}
if $u\in L^1(\Gw)\cap L^q_\gr(\Gw)$ and
\begin{equation}\label{int_form}
  -\int_\Gw u\Gd \vgf\,dx+ \int_\Gw |u|^q\sign u\vgf dx=-\int_\bdw\prtn\vgf\,d\nu,
\end{equation}
for every $\vgf\in C^2_0(\bar\Gw)$. If a solution exists it is unique; it will be denoted by $u_\nu$. If $\nu$ is a measure for which  a solution exists, we say that it is q-good. The family of q-good measures is denoted by $\CG_q$. It is known that \cite{BP84} $\nu$ is q-good if and only if it vanishes on sets of $\Capq$ capacity zero. Furthermore, a positive measure $\nu\in \GTM$ is q-good if and only if  it is the limit of an increasing
bounded \seq of measures in $\Wqq$. In particular a measure $\nu\in \GTM$ \sth $|\nu|\in \Wqq$ is a q-good measure.

A solution $u$ of \req{q-eq} is \emph{moderate} if $u\in L^1(\Gw)\cap L^q_\gr(\Gw)$. A moderate solution possesses a boundary trace $\nu\in \GTM$ \sth \req{int_form} holds.

Denote by $\CU_q$ the set of positive solutions of \req{q-eq}. A solution $u\in \CU_q$ is \emph{\gsmod} if it is the limit of an increasing \seq of moderate solutions.

A compact set $F\sbs \bdw$ is q-removable if a non-negative solution of \req{q-eq} vanishing on $\bdw\sms F$ must vanish in $\Gw$. An arbitrary set $A\sbs \bdw$ is q-removable if every compact subset is q-removable. It is known that $A$ is q-removable if and only if $\Capq(A)=0$
(see \cite{MVrem} and the references therein).

By \cite{MVcapest}, if $\{u_n\}$ is a \seq of positive solutions  of \req{q-eq} then
\begin{equation}\label{ac_Capq}
  \Capq(\suppq u_n)\to 0\Lra u_n\to 0 \text{ \locun in $\Gw$.}
\end{equation}

If $F$ is a \qcl subset of $\bdw$, denote
$$U_F=\sup\{u\in \CU_q:\,\suppq u\sbs F\}.$$
It is well known that $U_F$ is a solution of \req{q-eq} and it vanishes on $\bdw\sms F$. We call it \emph{the maximal solution relative to $F$.}

For an arbitrary Borel set $A\sbs \bdw$ denote
$$W_A=\sup\{ u_\nu: \, \nu\in \Wqq,\; \nu(\bdw\sms A)=0.\}$$
It is proved in \cite{MVcapest} that
$$W_A=W_{\tl A}$$
and, if $F$ is \qcl,
$$U_F=W_F.$$
In particular \emph{$U_F$ is \gsmod. }

If $v$ is a positive supersolution  of \req{q-eq} then the set of solutions dominated by it contains a maximal solution:
$$v^\#:=\sup\{u\in \CU_q:\, u\leq v\}\in \CU_q.$$

If $v$ is a positive subsolution  of \req{q-eq} then then the set of solutions dominating it is non-empty and contains a minimal solution:
$$v_\#:=\inf\{u\in \CU_q:\, u\geq v\}\in \CU_q.$$

If $u,v\in \CU_q$ then $u+v$ is a supersolution, $(u-v)_+$ is a subsolution  and we denote
$$u\oplus v=[u+v]^\#,\q u\ominus v=[(u-v)_+]_\#.$$

If $u\in \CU_q$ and $F$ is a \qcl subset of $\bdw$ we denote:
$$[u]_F=\inf(u,U_F)^\#.$$

If $D$ is a $C^2$ subdomain of $\Gw$ and $h\in L^1(\prt D)$ we denote by $S_q(D,h)$ the solution of the problem
$$-\Gd u+|u|^q\sign u=0 \text{ in $D$},  u=h\text{ on $\prt D$}.$$

Let $\{\Gw_n\}$ be a $C^2$ exhaustion of $\Gw$. Then, if $v$ is a positive supersolution,
$$S_q(\Gw_n,v)\;\dar \;v^\#$$
and if $v$ is a positive subsolution
$$S_q(\Gw_n,v)\;\uar \;v_\#.$$

The following definitions were introduced in \cite{MVfine}. A positive Borel measure $\tau$ on $\bdw$  (not necessarily bounded) is called a \emph{perfect measure} if it satisfies the following conditions:

\begin{quotation}
(a) $\tau$ is outer regular relative to $\Capq$-fine topology, i.e., for every Borel set $E$,
$$\tau(E)=\inf\{\tau(Q):\, Q \text{ is \qop},\;E\sbs Q\}.$$
(b) If $Q$ is a \qop set and $A$  a Borel set \sth $\Capq(A)=0$ then $\tau(Q)=\tau(Q\sms A)$.
\end{quotation}
The space of perfect measures is denoted by $\CM_q$.

We observe that (b) implies:
\begin{quotation}
(b') If $Q$ is a \qop set,  $A$ a Borel subset \sth $\Capq(A)=0$  and  $\tau (Q\sms A)<\infty$ then $\tau(A)=0$ and $\tau\chi_Q$ is a q-good measure.
\end{quotation}
For $\tau\in \CM_q$ put
$$Q_\tau=\bigcup\{Q:\,Q \text{ is \qop},\;\tau(\tl Q)<\infty\}.   $$
If $u\in \CU_q$ we say that $u$ has boundary trace $\tau\in \CM_q$
if:
\vskip 2mm
(i) $\CR(u)= Q_\tau$ and

(ii) for every $\gx\in Q_\tau$ there exists a \qop \ngh $Q$ \sth $[u]_{\tl Q}$ is a moderate solution with boundary trace $\tau\chi\ind{\tl Q}\,$.
\vskip 2mm
The boundary trace of $u$ in this sense is called \emph{the precise trace} and is denoted by $\tr u$ .

By \cite[Theorem 5.11]{MVfine}, for every $u\in \CU_q$, there exists a \seq $\{Q_n\}$ of \qop subsets of $\CR(u)$ \sth
$$\tl Q_n\sbs Q_{n+1},\q [u]_{\tl Q_n} \text{ is moderate} \forevery n,\q \Capq(\CR(u)\sms \cup_nQ_n)=0.$$
Such a \seq is called a \emph{regular decomposition} of $\CR(u)$. We denote:
\begin{equation}\label{CR_0}\BAL
 \CR_0(u)=\bigcup_n Q_n, &\q \nu_n=\tr [u]_{\tl Q_n},\\
u\ind{\CR}=\lim [u]_{\tl Q_n}, &\q \nu\ind{\CR}=\lim\nu_n.
\EAL
\end{equation}
$u\ind{\CR}$ and $\nu\ind{\CR}$ do not depend on the specific \seq $\{Q_n\}$. In fact (by the theorem cited above)
\begin{equation}\label{uR_F}
   [u]_F=[u\ind{\CR}]_F \forevery F\text{ \qcl,}\;F\sbs \CR(u),
\end{equation}
and $u\ominus u\ind{\CR}$ vanishes on $\CR(u)$.
\vskip 2mm

The following result  is proved in\cite{MVfine} (see Theorem 5.16 and the remark following it):

\bth{Ex+Un} Every positive solution $u$ of \req{q-eq} possesses a boundary trace $\nu\in \CB_q$. Conversely, for every $\nu\in \CB_q$ there exists a solution of \req{q-eq}  with boundary trace $\nu$. Furthermore there exists a unique \gsmod solution $u$ of \req{q-eq} with $\tr u=\nu$, namely,
$$u=u\ind\CR\oplus U_{\CS(u)}.$$
where $u\ind\CR$ is the \gsmod solution defined in \req{CR_0}.
\es

In addition, by \cite[Theorem 5.11]{MVfine} we obtain:

\bth{uRuS} If $u\in \CU_q$ then
\begin{equation}\label{uRuS}
\max(u_\CR,[u]\ind{\CS})\leq u\leq u\ind\CR +[u]\ind{\CS}.
\end{equation}
\es
\bproof $v:=u\ominus  u\ind\CR$ vanishes on $\CR(u)$, i.e., $\suppq v\sbs \CS(u)$. Thus $v$ is a solution dominated by $u$ and supported in $\CS(u)$, which implies that $v\leq [u]\ind\CS$. Since
$u- u\ind\CR\leq v$ this implies the inequality on the \rhs of \req{uRuS}. The inequality on the \lhs is obvious.
\eproof

We finish this section with the following lemma which is used in the proof of the main result.

\blemma{uAuB} Let $u\in \CU_q$  and let $A,B$ be two disjoint \qcl subsets of $\bdw$. If $u$ $\suppq u\sbs A\cup B$ and $[u]_A$, $[u]_B$ are \gsmod then $u$ is \gsmod. Furthermore
\begin{equation}\label{uAuB}
 u=[u]_A\oplus [u]_B=[\max(u_A,u_B)]_\#.
\end{equation}
\es
\bproof Let $\tau$ and $\tau'$ be q-good positive measures \sth $\qsupp \tau\cap\qsupp\tau'=\ems$. Then
$$[\max(u_\tau, u_{\tau'}]_\#= u_\tau\oplus u_{\tau'}=u_{\tau+\tau'}.$$
Let $\{\tau_n\}$ and $\{\tau'_n\}$ be increasing \seq{s} of q-good measures \sth
$$u_{\tau_n}\uar [u]_A, \q u_{\tau'_n}\uar [u]_B.$$
By \cite[Theorem 4.4]{MVfine}
\begin{equation}\label{max<u<sum}
   \max([u]_A,[u]_B)\leq u\leq [u]_A+[u]_B.
\end{equation}
Therefore
$$\max(u_{\tau_n}, u_{\tau'_n})\leq u\Lra  u_{\tau_n+\tau'_n}\leq u.$$
On the other hand
$$u-u_{\tau_n+\tau'_n}\leq ( [u]_A-u_{\tau_n})+([u]_B-u_{\tau'_n})\dar 0.$$
Thus
\begin{equation}\label{unun'}
 \lim u_{\tau_n+\tau'_n}= u
\end{equation}
so that $u$ is \gsmod.

Assertion \req{uAuB} is equivalent to  the statements:
(a) $u$ is the largest solution dominated by $[u]_A+[u]_B$ and (b) $u$ is the smallest solution dominating $\max(u_A,u_B)$.
Since the maximal solution $U_F$ of a \qcl set $F\sbs \bdw$ is \gsmod:
$$[u]_F=\sup\{v\in \CU_q:\, v\leq u,\; v \text{ moderate,}\; \suppq v\sbs F\}.$$

 Suppose that $w\in \CU_q$ and
$$u\leq w\leq [u]_A+[u]_B.$$
 Then,
$$[w]_A\leq [u]_A,\; [w]_B\leq [u]_B\Lra v \leq [u]_A.$$
Therefore, as $u\leq w$ we obtain,
$$[w]_A=[u]_A,\q [w]_B= [u]_B.$$
Since $u$ is \gsmod, these equalities and \req{unun'} imply that $u=w$.
This proves (a); statement (b) is proved in a similar way.
\eproof

\section{Characterization of positive solutions of $-\Gd u+u^q=0$.}

In this section we present the main result of the paper:

\bth{Main-q} Every positive solution of \req{q-eq} is \gsmod.
\es

The proof is based on several lemmas.

The following notation is used  throughout the section: $u$ is a positive solution of \req{q-eq},
$$V:=u^{q-1},\q L^V=-\Gd v+Vv=0.$$
Thus  $V$ satisfies \req{V-profile} and $L^Vu=0$. Therefore there exists a positive measure $\mu\in \GTM$ \sth
$$u=\BBK^V_\mu.$$
For any Borel set $E\sbs \bdw$ put
$$\mu_E=\mu\chi\ind E \;\text{ and }\; (u)_E=\BBK^V_{\mu_E}.$$

\blemma{mod1} Let $D$ be a $C^2$ domain \sth $D\Subset \Gw$ and let $h\in L^1(\prt D)$,
$0\leq h\leq u$. Then
\begin{equation}\label{mod1.1}
  S^V(D,h)\leq S_q(D,h).
\end{equation}
\es
\bproof Put $w:=S_q(D,h)$ and $v:=S^V(D,h)$. Then $w\leq u$ and \consy (recall that $V=u^{q-1}$)
$$0=-\Gd w+w^q\leq -\Gd w+Vw.$$
Thus $w$ is an  $L^V$ superharmonic  in $D$ \sth $u=h$ on $\prt D$. On the other hand $v$ is an $L^V$ harmonic in $D$ satisfying the same boundary condition. This implies \req{mod1.1}.
\eproof

\blemma{mod2} If $F$ is a compact subset of $\bdw$ then
\begin{equation}\label{(u)F}
   (u)_F\leq [u]_F.
\end{equation}
\es
\bproof Let $A$ be a Borel subset of $\bdw$. Put
$$A_\gb= \{x\in \Gw: \gr(x)=\gb,\;\gs(x)\in A\}$$
and
$$v_\gb^A=  S^V(D_\gb, u\chi\ind {A_\gb}), \q w_\gb^A=S_q(D_\gb,u\chi\ind {A_\gb}).$$
By \rlemma{mod1}, $v_\gb^A\leq w_\gb^A\leq u$. For any \seq $\{\gb_n\}$ decreasing to zero one can extract a \sseq $\{\gb_{n'}\}$ \sth
$\{w_{\gb_{n'}^A}\}$ and $\{v_{\gb_{n'}^A}\}$ converge locally uniformly; we denote the limits by $w^A$ and $v^A$ respectively. (The limits may depend on the \seq.) Then $w^A$ is a solution of \req{q-eq} while $v^A$ is an $L^V$ harmonic, and
\begin{equation}\label{vA<}
  v^A\leq w^A \leq  [u]_{\tl Q} \forevery Q \text{ open,}\; A\sbs Q.
\end{equation}
The second inequality follows from the fact that $w^A\leq u$ and $w^A$ vanishes on $\bdw\sms \tl Q$.

We apply the same procedure to the set $B=\bdw\sms A$ extracting a further \sseq of $\{\gb_{n'}\}$ in order to obtain the limits $v^B$ and $w^B$. Thus
$$v^B\leq w^B \leq  [u]_{\tl Q'} \forevery Q' \text{ open,}\; B\sbs Q'.$$

Note that
$$v^A+v^B=u,\q v^A\leq \BBK^V_{\mu_{\tl Q}},\q v^B\leq \BBK^V_{\mu_{\tl Q'}}.$$
Therefore
\begin{equation}\label{vA>}
  v^A=u-v^B\geq \BBK^V_{\mu\ind{\bdw\sms \tl Q'}}.
\end{equation}

Now, given $F$ compact, let $A$ be a closed set and $O$ an open set  \sth $F\sbs O\sbs A$ and let  $B=\bdw\sms A$. Note that $\bar B\cap F=\ems$. By \req{vA>} with $Q'=B$
$$v^A\geq  \BBK^V_{\mu\ind O}.$$
By \req{vA<}
$$v^A\leq w^A\leq [u]_{\tl Q} \forevery Q \text{ open,}\; A\sbs Q$$
and \consy
\begin{equation}\label{FOQ}
 (u)_F\leq \BBK^V_{\mu\chi\ind O}\leq  [u]_{\tl Q}.
\end{equation}
If $Q$ shrinks to $F$ then $[u]_{\tl Q}\dar [u]_F$ (see \cite[Theorem 4.4]{MVfine}).  Therefore \req{FOQ} implies \req{(u)F}.

\eproof
\blemma{q-mod1} If $E\sbs \bdw$ is a Borel set and $\Capq(E)=0$ then $\mu(E)=0$.
\es
\bproof If $F$ is a compact subset of $E$, $\Capq(F)=0$ and therefore the removability theorem \cite{MVrem} implies that $[u]_F=0$. Therefore, by \rlemma{mod2}, $(u)_F=0$. \Consy $\mu(F)=0$. As this holds for every compact subset of $E$ we conclude that $\mu(E)=0$.
\eproof

\blemma{q-mod2} Let $\nu\in W^{-2/q,q}(\bdw)$ be a positive measure and let $u_\nu$ be the solution of \req{q-eq} with trace $\nu$.  Suppose that there exists no positive solution of \req{q-eq} dominated by the supersolution $v=\inf(u,\BBK_\nu)$. Then $\mu\perp\nu$.
\es
\bproof
First we show:
\medskip

\note{Assertion 1} If
$V':=v^{q-1}$
then
$v$ is an $L^{V'}$ superharmonic and furthermore it is an $L^{V'}$ potential.
\medskip

Since $v$ is a supersolution of \req{q-eq}
$$0\leq -\Gd v+v^q= -\Gd v + V' v.$$
Thus $v$ is an $L^{V'}$ superharmonic. Suppose that there exists a positive $L^{V'}$ harmonic $w$ \sth $w\leq v$. Then $w$ is a subsolution of \req{q-eq}:
$$-\Gd w + w^q\leq -\Gd w + V' w=0.$$
This implies  that there exists a positive solution of \req{q-eq} dominated by $v$,
contrary to assumption. Thus $v$ is an $L^{V'}$-potential.

Note that
$$\int_{\bdw}\BBK_\nu V'\gr\, dx\leq \int_{\bdw}(\BBK_\nu )^q\gr\,dx<\infty.$$
Therefore
$\BBK_\nu$ is an  $L^{V'}$  superharmonic satisfying \req{LVr}. By  \rlemma{mod-basic} (i), the largest $L^{V'}$ harmonic dominated by $\BBK_\nu$, say $w$,  is $L^{V'}$ moderate and has m-boundary trace $\nu$. This implies that
$$ \BBK_\nu-w=:p$$
is an $L^{V'}$-potential. $w$ can be represented in the form
 $$w=\BBK^{V'}_{\nu'}$$ where $\nu'$ is a positive finite measure on $\bdw$ and, by \rprop{nu_nu'}, $\nu$, $\nu'$ are mutually a.c.

By the relative Fatou theorem, since $v,p$ are $L^{V'}$ potentials and $w$ is an $L^{V'}$ harmonic,
\begin{equation}\label{fine-nu}
 v/w\to 0, \q  \BBK_\nu/w\to1 \q L^{V'}-\text{finely $\nu'$-a.e.}
\end{equation}
Since $v=\inf(u,\BBK_\nu)$, \req{fine-nu} implies that
\begin{equation}\label{fine-nu'}
  u/w \to 0 \q \text{$L^{V'}$- finely $\nu'$-a.e.}
\end{equation}
Further, by  \req{fine-nu} and \req{fine-nu'}
\begin{equation}\label{fine-u}
 u/\BBK_\nu \to 0 \q \text{$L^{V'}$- finely $\nu'$-a.e.}
\end{equation}
Since $\nu$,$\nu'$ are mutually a.c.,  '$\nu$-a.e.' is equivalent to '$\nu'$-a.e.'. Therefore, in view of \rprop{fine-nt}, \req{fine-u} implies
\begin{equation}\label{nt-nu}
 u/\BBK_\nu \to 0 \q\text{n.t. $\nu$-a.e.}
\end{equation}

However, $\BBK_\nu$ is also an $L^V$ \superh. Therefore $\BBK_\nu$ can be represented in the form
$$\BBK_\nu=w^*+p^*,$$
where $w^*$ is an $L^V$-harmonic and $p^*$ an $L^V$-potential.
Let $\tau\in \GTM$ be the $L^V$ trace of $w^*$, i.e.,  $w^*=\BBK^V_\tau$. Then, by the relative Fatou theorem,
$$\BBK_\nu/u\to \frac{d\tau}{d\mu}=: h\q L^V-\text{finely,\;$\mu$-a.e.}$$
and therefore, by \rprop{fine-nt},
\begin{equation}\label{nt-tau}
  \BBK_\nu/u\to h\q\text{n.t. $\mu$-a.e.}
\end{equation}
Since $0\leq h<\infty$ $\mu$-a.e., \req{nt-nu} and \req{nt-tau} imply that $\nu\perp\mu$.
\eproof

\blemma{q-mod3}   Suppose that for every positive measure $\nu\in W^{-2/q,q}(\bdw)$, there exists no positive solution of \req{q-eq} dominated by $v=\inf(u,\BBK_\nu)$. Then $u=0$.
\es
\bproof By \rlemma{q-mod2},
$$\mu\perp\nu \forevery \nu\in W^{-2/q,q}(\bdw), \;\nu\geq 0.$$
Suppose that $\mu\neq 0$. By \rlemma{q-mod1}, $\mu$ vanishes on sets of $\Capq$ zero. Therefore (by Feyel and de la Pradelle \cite{FdlP} or Dal Maso \cite{DMaso}) $\mu$ is the limit of an increasing \seq $(\mu_k)\sbs \Wqq_+$.
For every $k$ there exists a Borel set  $A_k\sbs \bdw$ \sth,
$$\mu(A_k)=0,\q \mu_k(\bdw\sms A_k)=0.$$
Therefore, if $A=\cup A_{k}$ and $A'=\bdw\sms A$ then
$$\mu(A)=0, \q\mu_k(A')=0 \forevery k.$$

Since $\mu_k\leq\mu$ we have $\mu_k(A)=0$ and therefore $\mu_k=0$. Contradiction!
\eproof

\vskip 2mm

\noindent\textsc{Proof of \rth{Main-q}}  Let $\{Q_n\}$ be a regular decomposition of $\CR(u)$ and put
$$v_n:=[u]_{\tl Q_n}.$$
Using the notation introduced in \req{CR_0}, $v_n$ is moderate with boundary trace $\nu_n$ and
$$v_n\uar u_\CR.$$
Thus the solution $u_\CR$ is \gsmod and
$$u\ominus u_\CR \leq [u]_{\CS(u)}=:u_\CS.$$
\vskip 2mm

\emph{Assertion 1} $u_\CS$ is \gsmod.
\vskip 2mm

 Before proving the assertion we verify that it implies that $u$ is \gsmod. Put
 $$u_n:=v_n\oplus u_\CS.$$
By \rlemma{uAuB}, as $\tl Q_n\cap \CS(u)=\ems$, it follows that $u_n$ is \gsmod.
As $\{u_n\}$ is increasing it follows that $\bar u=\lim u_n$ is a \gsmod solution of \req{q-eq}.
In addition
$$[\max(v_n, u_\CS)]_\#=u_n=v_n\oplus u_\CS \Lra \max(u_\CR, u_\CS)\leq \bar u\leq u_\CR+u_\CS.$$
This further implies that $\CS(u)=\CS(\bar u)$ and that $\tr \bar u=\tr u$. By uniqueness of the \gsmod solution we conclude that $u=\bar u$.
\vskip 2mm

We turn to the proof of \emph{Assertion 1}.  To simplify notation, we put  $u=u_\CS$  and denote $F:=\suppq u$. (Incidentally, $F\sbs \CS(u)$ but it is possible that there is no equality. In fact $F$ consists precisely of the $\Capq$-thick points of $\CS(u)$. The set $\CS(u)\sms F$ is contained in the singular set of $u_\CR$.)

For $\nu\in \Wqq$ we denote by $u_\nu$ the solution of \req{q-eq} with boundary trace $\nu$. Put
\begin{equation}\label{u*}
   u^*=\sup\{u_\nu:\, \nu\in \Wqq,\; 0<u_\nu\leq u\}.
\end{equation}
By \rlemma{q-mod3} the family over which the supremum is taken is not empty. Therefore $u^*$ is a positive solution of \req{q-eq} and it is well-known that it is \gsmod. By its definition, $u^*\leq u$. 

Let $F^*=\suppq u^*$. Then $F^*$ is \qcl and $F^*\sbs F$. Suppose that $\Capq(F\sms F^*)>0$. Then there  exists a compact set $E\sbs F\sms F^*$ \sth $\Capq(E)>0$ and $\bdw\sms F^*=:Q^*$ is a \qop set containing $E$. Furthermore there exists a \qop set $Q'$ \sth $E\sbs Q'\sbs \tl Q'\sbs Q^*$ (\cite[Lemma 2.4]{MVfine}) .
Since $Q'\sbs \suppq u$,  $[u]_{\tl Q'}>0$ and therefore, by \rlemma{q-mod3}, there exists a positive measure $\tau\in \Wqq$ supported in $\tl Q'$  \sth
$u_\tau\leq u$. As the $\qsupp\tau$ is a \qcl set disjoint from $F^*$ it follows that $u^*\ngeq u_\tau$. On the other hand, since $\tau\in \Wqq$ and $u_\tau\leq u$, it follows that $u_\tau\leq u^*$. This contradiction shows that
\begin{equation}\label{FF*}
  \Capq(F\sms F^*)=0.
\end{equation}

Further $u^*$ is \gsmod and therefore there exists a \qcl set $F_0^*\sbs F^*$ \sth $\CS(u^*)=F_0^*$ and
$\CR(u^*)=\bdw\sms F_0^*$.
Suppose that $\Capq(F\sms F_0^*)>0$ and put $Q_0:=\bdw\sms F_0^*$. Let $E\sbs F\sms F_0^*$ be a compact set \sth $\Capq(E)>0$ and let $Q'$ be a \qop set \sth $E\sbs Q'\sbs \tl Q'\sbs Q_0$. Then $\tl Q'\sbs \CR(u^*)$ and \consy $[u^*]_{\tl Q'}$ is a moderate solution of \req{q-eq}, i.e.
$$[u^*]_{\tl Q'}\in L^q_\gr(\Gw).$$
On the other hand $Q'$ is a \qop \ngh of $E$ which is a non-empty subset of $F=\suppq u$; therefore $[u]_{\tl Q'}$ is a purely singular solution of \req{q-eq}, i.e.,
$$\int_\Gw ([u]_{\tl Q'})^q\gr\,dx=\infty,\q \CS([u]_{\tl Q'})=\suppq [u]_{\tl Q'}.$$
It follows that
$v:=\big[[u]_{\tl Q'}-[u^*]_{\tl Q'}\big]_\#$ is a purely singular solution  of \req{q-eq}.

Let $v^*$ be defined as in \req{u*} with $u$ replace by $v$.  Then $v^*$ is a singular, \gsmod solution of \req{q-eq}. Since $v^*\leq u$ and it is \gsmod it follows that $v^*\leq u^*$. On the other hand, since $v^*$ is singular and $\suppq v^*\sbs \tl Q'\sbs \CR(u^*)$ it follows that $u^*\ngeq v^*$, i.e. $(v^*-u^*)_+$ is not identically zero. Since both $u^*$ and $v^*$ are \gsmod, it follows that there exists $\tau\in \Wqq$ \sth $u_\tau\leq v^*$ but $(u_\tau- u^*)_+$ is not identically zero. Therefore $u^*\lneq\max(u^*,u_\tau)$. The function $\max(u^*,u_\tau)$ is a subsolution of \req{q-eq} and the smallest solution above it, which we denote by $Z$ is \emph{strictly larger} then $u^*$. However  $u_\tau\leq v^*\leq u^*$ and \consy $Z=u^*$. This contradiction proves that
\begin{equation}\label{FF*0}
  \Capq(F\sms F_0^*)=0
\end{equation}

In conclusion,  $u^*$ is \gsmod, $\suppq u^*\sbs F$ and $F_0^*=\CS(u^*)$  is $\Capq$-equivalent to $F$.Therefore, by \rth{Ex+Un},$u^*=U_F$, the maximal solution supported in $F$. Since, by definition $u^*\leq u$, it follows $u^*=u$.
\qed


\begin{thebibliography}{99}

\bibitem{AH} Adams D. R. and Hedberg L. I., Function spaces and potential theory,
Grundlehren  Math. Wissen. {\bf 314}, Springer (1996).
\bibitem{An-SLN} Ancona A., \emph{Theorie du potentiel sur les graphes et les varietes,}
in Springer Lecture Notes No. 1427 (ed. P.L. Hennequin) p.1-112 (1988).
\bibitem{An87} Ancona A., \emph{Negatively curved manifolds, elliptic operators and the Martin boundary,} Annals of Mathematics, Second Series, {\bf 125}, 495-536 (1987).
\bibitem{An-05} Ancona A., Personal communication, 2005.
\bibitem{An-10} Ancona A., Personal communication, 2010.
\bibitem{An-App} Ancona A., \emph{A necessary condition for the fine regularity of a boundary point
with respect to a Schr\"odinger equation,} Appendix in `Boundary value problems with measures', by L.V\'eron and C.Yarur (preprint).
\bibitem{BM90} Bandle C. and Marcus M. \emph{Sur les solutions maximales de problèmes elliptiques non linéaires,} C. R. Acad. Sci. Paris Ser. I {\bf 311} (1990) 91–93.
\bibitem{BM95} Bandle C. and Marcus M. \emph{Asymptotic behaviour of solutions and their derivatives for semilinear
elliptic problems with blow up on the boundary,} Ann. Inst. Poincaré {\bf 12} (1995) 155–171.
\bibitem{BP84} Baras P., Pierre M. \emph{Singularit\`es \'eliminables pour des \'equations semilin\`eaires}, Ann. Inst. Fourier 34 (1984), 185--206.
\bibitem{DMaso} Dal Maso G. \emph{On the integral representation of certain local functionals,} Ricerche Mat. {\bf 32 no. 1} (1983)  85-113.
\bibitem{Dy91} Dynkin E. B., \emph{A probabilistic approach to one class of nonlinear differential equations}, Probab. Th, Rel. Fields 89 (1991), 89--115.
\bibitem{Dy93} Dynkin E. B., \emph{Superprocesses and partial differential equations}, Ann. Probab 21 (1993), 1185--1262.
\bibitem{Dy98} Dynkin E. B., \emph{Stochastic boundary values and boundary singularities for solutions of the equation $Lu = u^\alpha$}, J. Funct. Anal. 153 (1998), 147--186.
\bibitem{Dbook1} Dynkin E. B. {\em Diffusions, Superdiffusions and Partial Differential Equations},
American Math. Soc., Providence, Rhode Island, Colloquium Publications {\bf 50}, 2002.
\bibitem{Dbook2}  Dynkin E. B. {\em Superdiffusions and Positive Solutions of Nonlinear
Partial Differential Equations},
American Math. Soc., Providence, Rhode Island, Colloquium Publications {\bf 34}, 2004.
\bibitem {DK96} Dynkin E. B. and Kuznetsov S. E. {\em Superdiffusions
and removable singularities for quasilinear partial differential
equations}, Comm. Pure Appl. Math. {\bf 49}, 125-176 (1996).
\bibitem {DK98a} Dynkin E. B. and Kuznetsov S. E. {\em Trace on the boundary for solutions of nonlinear differential
equations,} Trans. Amer. Math. Soc. {\bf 350} (1998) 4499-4519.
\bibitem{DK98b} Dynkin E. B. and Kuznetsov S. E. {\em Fine topology and fine trace on the boundary associated with
a class of quasilinear differential equations}, Comm. Pure Appl. Math. {\bf 51}, 897-936 (1998).
\bibitem{FdlP}  Feyel D. and de la Pradelle A. \emph{Topologies fines et compactifications associes certains
espaces de Dirichlet,} Ann. Inst. Fourier (Grenoble) 27, 121146 (1977).
\bibitem{GV} Gmira A. and V´eron L.,\emph{Boundary singularities of solutions of nonlinear elliptic
equations,} Duke J. Math. 64, 271-324 (1991).
\bibitem{Kz98} Kuznetsov S.E. {\em \gsmod solutions of $Lu=u^\ga$ and fine trace on the boundary}, C.R.Acad. Sc. Serie I {\bf 326}, 1189-1194 (1998).
 \bibitem{Kell} Keller   J. B. , \emph{On solutions of $\Delta u = f(u)$}, Comm. Pure Appl. Math. 10, 503-510
(1957).
\bibitem{LG93a} Le Gall J. F., {\em A class of path-valued Markov processes and its connections with partial
differential equations.,} Probab. Th. Rel. Fields 102 (1993), 25–46.
\bibitem{LG93b} Le Gall J. F., {\em Solutions positives de $\Delta u = u^2$ dans le disque unit\`e,} C.R. Acad. Sci. Paris, S\'er. I 317 (1993), 873–878.
\bibitem{LG95} Le Gall J. F., {\em The Brownian snake and solutions of
$\Gd u=u^{2}$ in a domain}, Probab. Th. Rel. Fields {\bf 102}, 393-432 (1995).
\bibitem{LG97} Le Gall J. F., {\em A probabilistic Poisson representation for positive solutions of $\Delta u = u^2$ in a
domain,} Comm. Pure Appl. Math. 50 (1997), 69–103.
\bibitem{LGbook} Le Gall J. F., {\em Spatial branching processes, random snakes and partial differential equations}, Birkh\"{a}user, Basel/Boston/Berlin, 1999.
\bibitem{LN} Loewner C. and Nirenberg L. \emph{Partial differential equations invariant under conformal or
projective transformations,} Contributions to Analysis, Academic Press, Orlando, FL (1974),
245–272.
\bibitem{MVunique} Marcus M. and V\'{e}ron L., {\em Uniqueness and asymptotic behaviour of solutions with
boundary blow-up for a class of nonlinear elliptic equations,} Annales de l'Institut
Henri Poincard 14 (1997), 237-274.
\bibitem{MVcras96} Marcus M. and V\'{e}ron L., {\em Trace au bord des solutions positives d'\'equations elliptiques et
paraboliques non lin\'eaires. R\'esultats d'existence et d'unicit\'e,}  [Boundary trace of positive
solutions of nonlinear parabolic and elliptic equations. Existence and uniqueness results.] C.
R. Acad. Sci. Paris S\'er. I, Math. {\bf 323, no.6} (1996) 603–608.
\bibitem{MVsub} Marcus M. and V\'{e}ron L., {\em The boundary trace of positive
solutions of semilinear elliptic equations: the subcritical
case}, Arch. rat. Mech. Anal. {\bf 144}, 201-231 (1998).
\bibitem{MVsuper} Marcus M. and V\'{e}ron L., {\em The boundary trace of positive
solutions of semilinear elliptic equations: the supercritical
case}, J. Math. Pures Appl. {\bf 77}, 481-524 (1998).
\bibitem{MVrem} Marcus M. and V\'{e}ron L., {\em Removable
singularities and boundary trace}, J. Math. Pures Appl. {\bf 80},
879-900 (2000).
\bibitem{MVcapest} Marcus M. and V\'{e}ron L., {\em Capacitary estimates of positive solutions of semilinear elliptic equations with absorption}, J. European Math. Soc. {\bf 6}, 483-527 (2004).
\bibitem{MVfine} Marcus M. and V\'{e}ron L., \emph{The precise boundary trace of positive solutions of the equation $\Gd u=u^q$ in the supercritical case,} in "Perspectives in Nonlinear Partial Differential
Equations" (in honor of H. Brezis), Berestycki et al eds., Contemporary
Mathematics  A. M. S. {\bf 446} (2007) p. 345-384.
\bibitem{MVbesov} Marcus M. and  V´eron L., \emph{On a  characterisation of Besov spaces with negative exponents,}
Topics Around the Research of Vladimir Mazya /International Mathematical
Series 11-13, Springer Science+Business Media, 2009.
\bibitem{MVcapint} Marcus M. and  V´eron L., \emph{Maximal solutions for $-\Delta u+u^q=0$ in open and finely open sets,} J. Math. Pures Appl. {\bf 91} (2009) 256-295.
\bibitem{Ms} Mselati B., \emph{Classification and probabilistic representation of the positive solutions of a semilinear elliptic equation}, (English summary)
Mem. Amer. Math. Soc. {\bf 168, no. 798},  (2004) 121 pp.
\bibitem{Oss} Osserman R., \emph{On the inequality $\Delta u\geq f(u)$,} Pacific J. Math. 7, 1641-1647 (1957).
    \bibitem{Triebel} Triebel H., Interpolation Theory, Function Spaces, Differential Operators, North-Holland Pub. Co., 1978.
\end{thebibliography}
\end{document}